\documentclass[12pt]{amsart}
\usepackage{amssymb, amsmath, amsfonts}
\usepackage[raiselinks=false,colorlinks=true,citecolor=blue,urlcolor=blue,linkcolor=blue,bookmarksopen=true,pdftex]{hyperref}
\newtheorem{theorem}{Theorem}[section]
\newtheorem{proposition}[theorem]{Proposition}
\newtheorem{lemma}[theorem]{Lemma}
\newtheorem{remark}[theorem]{Remark}

\newcommand{\out}{\textrm{Out}}

\newcommand{\bz}{\mathbb{Z}}
\newcommand{\bq}{\mathbb{Q}}
\newcommand{\br}{\mathbb{R}}
\newcommand{\bc}{\mathbb{C}}

\newcommand{\lr}{\longrightarrow}

\newcommand{\wt}{\widetilde}

\newcommand{\Sp}{\textrm{Sp}}

\newcommand{\gl}{\textrm{GL}}
\newcommand{\SL}{\textrm{SL}}
\newcommand{\pgl}{\textrm{PGL}}
\newcommand{\psl}{\textrm{PSL}}
\newcommand{\tr}{\textrm{tr}}

\begin{document}
\baselineskip=15.5pt
\title[Twisted conjugacy classes]{Twisted conjugacy classes in abelian extensions of certain linear 
groups}  
\author{T. Mubeena}
\author{P. Sankaran}
\address{The Institute of Mathematical Sciences, CIT
Campus, Taramani, Chennai 600113, India}
\email{mubeena@imsc.res.in}
\email{sankaran@imsc.res.in}


\date{}

\begin{abstract}
 Given an automorphism $\phi:\Gamma\lr \Gamma$,  one has 
an action of $\Gamma$ on itself by $\phi$-twisted conjugacy, namely, $g.x=gx\phi(g^{-1})$. 
The orbits of this action are called $\phi$-twisted conjugacy classes.  One says 
that $\Gamma$ has the $R_\infty$-property if there are infinitely many $\phi$-twisted conjugacy 
classes for every automorphism $\phi$ of $\Gamma$. In this paper we show that $\SL(n,\bz)$ and its 
congruence subgroups have the $R_\infty$-property.  Further we show that 
any (countable) abelian extension of $\Gamma$ has the $R_\infty$-property where $\Gamma$ is a torsion 
free non-elementary hyperbolic group, or $\SL(n,\bz), \Sp(2n,\bz)$ or a principal congruence 
subgroup of $\SL(n,\bz)$ or the fundamental group of a complete Riemannian 
manifold of constant negative curvature. 
\end{abstract}
\maketitle
\section{Introduction}
Let $\Gamma$ be a finitely generated infinite group and let $\phi:\Gamma\lr \Gamma$ be an endomorphism.  
One has an action of $\Gamma$ on itself defined as $g.x=gx\phi(g^{-1})$.  This is just the conjugation action 
when $\phi$ is identity. The orbits of this action are called the $\phi$-twisted conjugacy classes; the $\phi$-twisted conjugacy class containing $x\in \Gamma$ is denoted $[x]_\phi$ or simply $[x]$ when $\phi$ is clear from the context. If $x$ and $y$ are in the same $\phi$-twisted conjugacy class, we write $x\sim_\phi y$.  The set of all $\phi$-twisted conjugacy classes is denoted by $\mathcal{R}(\phi)$. The cardinality $R(\phi)$ of $\mathcal{R}(\phi)$ is called the Reidemeister number of 
$\phi$.  One says that $\Gamma$ has the $R_\infty$-property for automorphisms (more briefly, $R_\infty$-property) if there are infinitely many $\phi$-twisted conjugacy classes for every automorphism $\phi$ of $\Gamma$.  If $\Gamma$ has the $R_\infty$-property, we shall call $\Gamma$ an $R_\infty$-group.  All these notions make sense for any group, not necessarily finitely generated. 

The notion of twisted conjugacy originated in 
Nielson-Reidemeister fixed point theory and also arises in other areas of mathematics such as representation theory, number theory and algebraic geometry. See \cite{felshtyn2} and the references therein. 
The problem of determining which classes of groups have $R_\infty$-property is an area of active research initiated by Fel'shtyn and Hill \cite{fh}.  
We now state the main result of this paper.     
 
\begin{theorem} \label{main}
Let $\Lambda$ be an extension of a group $\Gamma$ by an arbitrary countable abelian group $A$. 
Then $\Lambda$ has the $R_\infty$-property in case any one of the following holds:\\
(i) $\Gamma$ is a torsion-free non-elementary hyperbolic group,\\
(ii) $\Gamma=\SL(n,\bz), \psl(n,\bz), \pgl(n,\bz), \Sp(2n,\bz),\textrm{PSp}(2n,\bz)$, $n\geq 2,$\\ 
(iii) $\Gamma$ is a normal subgroup of $\SL(n,\bz), n>2,$ not contained in the centre, and, \\  
(iv) $\Gamma$ is the fundamental group of a complete Riemannian manifold of constant negative sectional curvature and finite volume. 
\end{theorem} 

Our proofs involve straightforward arguments, using well-known results concerning the group $\Gamma$ in each case. More precisely,  in each of the cases, we show that $A$ or a bigger subgroup $N\subset \Lambda$ in which $A$ 
has finite index is {\it characteristic} in $\Lambda$. Proof of 
this requires some facts concerning normal subgroups 
of $\Gamma$. In the cases (ii), (iii) and (iv) we invoke 
the normal subgroup theorem of Margulis \cite[Chapter 8]{zimmer}; in case (i) we use the quasi-convexity property 
of infinite cyclic subgroups of $\Gamma$. 
Using the fact that $\Gamma$ is hopfian, the $R_\infty$-property for 
$\Lambda$ is then deduced from the $R_\infty$-property for 
$\Gamma$.  That $\Gamma$ has the $R_\infty$-property when it is a torsion-free non-elementary hyperbolic group is due to \cite{ll}. This  result was extended to arbitrary non-elementary hyperbolic groups by \cite{felshtyn1}.  The $R_\infty$-property for the groups $\Sp(2n,\bz)$ was established by Fel'shtyn and Gon\c{c}alves \cite{fg}.    
The $R_\infty$-property for $\SL(n,\bz)$ and $\pgl(n,\bz)$ is established in \S3.  We show, in \S3, the $R_\infty$-property for non-central normal subgroups of $\SL(n,\bz), n>2,$ using the Mostow-Margulis strong rigidity theorem and the congruence subgroup property of $\SL(n,\bz)$. The proof of the main theorem is given in \S4.  
\section{preliminaries}

Let $G$ be a group and $H$ a subgroup of $G$.  Recall that a subgroup $H$ is said to be {\it  characteristic} in $G$ if $\phi(H)=H$ for every automorphism $\phi$ of $G$.  $G$ is called {\it hopfian} (resp. {\it co-hopfian}) if every surjective (resp. injective) endomorphism of $G$ is an automorphism of $G$.   One says that $G$ is {\it residually finite}  
if, given any $g\in G$, there exists a finite index subgroup $H$ in $G$ such that $g\notin H$.

We shall recall here some facts concerning the $R_\infty$-property. 
Let \[1\lr N\stackrel{j}{\hookrightarrow} \Lambda\stackrel{\eta}{\lr}\Gamma\lr 1\eqno(1)\] be an exact sequence of groups.    
  
\begin{lemma} \label{characteristic} 
Suppose that $N$ is characteristic in $\Lambda$ and that $\Gamma$ has the $R_\infty$-property, then $\Lambda$ also has the $R_\infty$-property.
\end{lemma}

\begin{proof}
Let $\phi:\Lambda\lr \Lambda$ be any automorphism. Since $N$ is characteristic, $\phi(N)=N$ and so $\phi$ induces an automorphism $\bar{\phi}:\Gamma\lr \Gamma$. Since $R(\bar{\phi})=\infty$, it follows that $R(\phi)=\infty$.   
\end{proof}

\begin{lemma}\label{finitecharacteristic}
Suppose that $N$ is a characteristic subgroup $\Lambda$. 
(i) If $N$ is finite and $\Lambda$ has the $R_\infty$-property then $\Gamma$ also has the $R_\infty$ 
-property.\\
(ii) If $\Gamma$ is finite and $N$ has the $R_\infty$-proeperty, then $\Lambda$ has the $R_\infty$- property.  
\end{lemma}

\begin{proof}
(i)  Any automorphism $\phi:\Lambda\lr \Lambda$ maps $N$ isomorphically onto itself and hence induces an 
automorphism $\bar{\phi}:\Gamma \lr \Gamma$ (where $\Gamma=\Lambda/N$).

It is readily seen that $x\sim_\phi y$ implies $\eta(x)\sim_{\bar{\phi}}\eta(y)$ for any $x,y\in \Lambda$. 
Therefore $\eta$ induces a surjection  
$\wt{\eta}:\mathcal{R}(\phi)\lr \mathcal{R}(\bar{\phi})$ where $\wt{\eta}([x]_\phi)=[\eta(x)]_{\bar{\phi}}$. 
We need only show that the fibres of $\wt{\eta}$ are 
finite. 

Suppose the contrary and let $x_k\in \Lambda, k\geq 0,$ be such that $[x_k]_\phi\neq [x_l]_\phi$ for $k\neq l$ and that $[\eta(x_k)]_{\bar{\phi}}= [\eta(x_0)]_{\bar{\phi}}$.  For each $k\geq 1$, there exists $g_{k}\in \Lambda$ such that 
$\eta(x_0)=\eta(g_{k})\eta(x_k)\bar{ \phi}(\eta(g_{k}^{-1}))
=\eta(g_kx_k\phi(g_k^{-1}))$.  Therefore there exists an $h_k\in N$ such that 
$x_0h_k=g_{k}x_k\phi(g_{k})^{-1}$.  That is, for any $k\geq 1$, we have $x_k\sim_\phi x_0h_k$ for some $h_k\in N$. Since $N$ is finite, it follows that $x_k\sim_\phi x_l$ for some $k\neq l$, a contradiction. 

(ii)  Let $\phi:\Lambda\lr \Lambda$ be an automorphism and let $\theta=\phi|N$.  Let $\wt{j}:\mathcal{R}(\theta)\lr \mathcal{R}(\phi)$ be the map defined as 
$[x]_\theta\mapsto [x]_\phi$.  Suppose that $R(\phi)<\infty$
but that $R(\theta)=\infty.$  Then there exist elements $x_k\in N, k\geq 0,$ such that $[x_k]_\theta\neq [x_l]_\theta, k\neq l,$ but  
$x_k\sim_\phi x_0$ for all $k\geq 0$.   Choose $g_k\in \Lambda$ such that $x_k=g_kx_0\phi(g_k^{-1}), k\geq 1$.  
Since $\Gamma=\Lambda/N$ is finite, there exist distinct  positive integers $k,l$ such that $h:=g_kg_l^{-1}\in N$. 
Now $x_k=g_kx_0\phi(g_k^{-1})=g_kg_l^{-1}x_l\phi(g_l)\phi(g_k^{-1})=hx_l\theta(h^{-1})$ and so $[x_k]_\theta
=[x_l]_\theta$, a contradiction. This completes the proof.
\end{proof}   

\begin{lemma} \label{hopfian}
Suppose that there is no non-trivial  homomorphism from $N$ to $\Gamma$ and that either $\Gamma$ is hopfian or $N$ is co-hopfian. If $\Gamma$ has the $R_\infty$-property, then so does $\Lambda$. 
\end{lemma}
\begin{proof}
Let $\phi:\Lambda\lr \Lambda$ be any automorphism. Consider the homomorphism $f:N\lr \Gamma$ defined as $f=\eta\circ \phi|N$ where $\eta:\Lambda\lr \Gamma$ is the quotient map as in (1). By our hypothesis $f$ is trivial, and so it follows that $\phi(N)\subset \ker(\eta)=N$.  If $N$ is co-hopfian then $\phi(N)=N$ and so $N$ is characteristic. In any case $\phi$ defines a homomorphism $\bar{\phi}:\Gamma\lr \Gamma$ where $\bar{\phi}(xN)=\phi(x)N, ~x\in \Lambda$. It is clear that $\bar{\phi}$ is surjective with kernel $\phi^{-1}(N)/N$. If $\Gamma$ is hopfian, $\bar{\phi}$ is an isomorphism and it follows that 
$\phi(N)=N$. Thus our hypothesis implies that  $N$ is characteristic in $\Lambda$ and the lemma now follows from Lemma \ref{characteristic}. 
\end{proof}

We conclude this section with the following observation.

\begin{proposition}\label{residual}
Let $\Gamma$ be a countably infinite residually finite group. Then $R(\phi)=\infty$ for any inner automorphism $\phi$ of $\Gamma$.
\end{proposition}

\begin{proof}
Let $\phi=\iota_\gamma$ and let $x\sim_\phi y.$  Thus $y=gx\gamma g^{-1}\gamma^{-1}$.  Equivalently $x\gamma$ is conjugate to $y\gamma$. 
Hence it suffices to show that $\Gamma$ has infinitely many conjugacy classes.  

Since $\Gamma$ is infinite and since $\Gamma$ is residually finite, 
there exist finite quotients $\bar{\Gamma}$ of $\Gamma$ having {\it arbitrarily large} (finite) order. 
It is a classical result of R. Brauer \cite{brauer} (see also \cite{pyber})
that the number of conjugacy classes of a finite group of order $n$ is bounded below by $\log \log n$.  Since $\Gamma$ has at least as many conjugacy classes as any of its quotients, it follows that $\Gamma$ has infinitely many conjugacy classes.  
\end{proof}  

\begin{remark}{\em (i) Recall that finitely generated 
residually finite groups are hopfian.  A well-known 
class of residually finite groups is the class of finitely 
generated subgroups of $\gl(n, K)$ where $K$ is any field. See \cite{ls}. This class includes, in particular,  all lattices in linear Lie groups.  An important unsolved problem is to decide whether hyperbolic groups are 
residually finite. It has been shown by Sela \cite{sela} that 
torsion-free hyperbolic groups are hopfian. 
 
(ii) It is known that there are countably infinite groups with only finitely many conjugacy classes. (See \cite[\S 1.4]{serre2}, or  \cite[Chapter 4, \S 3]{ls}.) Finitely generated such examples have been constructed by S. Ivanov.  Recently D. Osin \cite{osin} has constructed a finitely generated infinite group which has exactly one non-trivial conjugacy class.}
\end{remark}
\section{The $R_\infty$-property for special linear groups}
In this section we give new examples of $R_\infty$-groups which are subgroups of finite index in 
$\SL(n,\bz)$. 
Recall that $\Sp(2n;\bz)$ has been shown to have the $R_\infty$-property by Fel'shtyn and Gon\c{c}alves \cite{fg}.

Our first result is the following.

\begin{theorem}\label{linear}
The groups $\SL(n,\bz), \psl(n,\bz), \gl(n,\bz)$, and $\pgl(n,\bz)$ have the $R_\infty$-property for all $n\geq 2$.   
\end{theorem}
\begin{proof}
It follows from Lemma \ref{finitecharacteristic}(ii) that the $R_\infty$-property for $\SL(n,\bz)$ implies 
that for $\gl(n,\bz)$.  Also the $R_\infty$-property 
for $\SL(n,\bz)$ (resp. $\gl(n,\bz)$) implies that for $\psl(n,\bz)$ (resp. $\pgl(n,\bz)$) in view of Lemma \ref{characteristic}.   Therefore we need 
only prove the theorem for $\SL(n,\bz).$

The group $\SL(2,\bz)$ is non-elementary hyperbolic group and hence, by \cite{felshtyn1}, has the $R_\infty$-property.  
Let $n\geq 3$ and set $\Gamma:=\SL(n,\bz)$.  
Since $\Gamma$ is residually finite, $R(\phi)=\infty$ for any inner automorphism by Proposition \ref{residual}.  In this case we can see this more directly: the set $ \{tr(A)\mid A\in \Gamma\}$ is infinite and so there are infinitely many conjugacy classes in $\Gamma$. 

It remains only to show that $\mathcal{R}(\phi)$ is infinite for a set of representatives of (the non-trivial) elements of the group $\out(\Gamma)$ of all outer automorphisms of $\Gamma$. 
It is known from the work of Hua and Reiner \cite{hr} and of O'Meara \cite{omeara} that $\out(\Gamma)\cong \bz/2\bz$ or $\bz/2\bz \times \bz/2\bz$ according as $n$ is odd or even.  

The group $\out(\Gamma)$ is generated by a set $S$ where $S=\{\tau\}$ when $n$ is odd and when $n$ is even, $S=\{\sigma, \tau\}$ where $\tau:\Gamma\lr \Gamma$ is defined as $X\mapsto {}^tX^{-1}$, and, when $n$ is even, the 
involution $\sigma:\Gamma\lr \Gamma$ is defined as 
$X\mapsto JXJ^{-1}=JXJ$ where $J$ is the diagonal matrix 
$diag(1,\ldots, 1,-1)$.  Thus $X\sim_\tau Y$ (resp. $X\sim_\sigma Y$) if and only 
if there exists a $Z$ such that $Y= ZX({}^tZ)$ (resp. $Y=ZXJZ^{-1}J)$).

First we consider $\tau$-twisted conjugacy classes.  
Let $k\geq 1$ and let $A(k)$ be the 
block diagonal matrix $A(k)=diag(B(k), I_{n-2})$ where $B(k)=\left(\begin{array}{cc} 1 & 0\\k&1\\
 \end{array}\right)$.  We shall show that $A(k)\sim_\tau A(l)$ implies $k=l$.  This clearly implies that $R(\tau)=\infty$.  
 
Let $X=(x_{ij})\in \Gamma$ be such that
               \[X.A(k).^tX=A(l).    \eqno(2)\]  
We shall denote the $i$-th row and $i$-th column of $X$ by $r_i$ and $c_i$ respectively.    
A straightforward computation shows that 
$X.A(k).^tX = X.^tX+kc_2.{}^tc_1$. 
Comparing the $(2,1)$-entries on both sides of (2) we get $r_2.{}^tr_1+ kx_{22}x_{11}=l$ whereas comparing the $(1,2)$-entries gives $r_1.{}^tr_2+k x_{12}x_{21}=0$.  Therefore $r_2.{}^tr_1=r_1.{}^tr_2= -kx_{12}x_{21}$ and so $l=k(x_{11}x_{22}-x_{12}x_{21})$.  Since $x_{i,j}\in \bz$, we obtain that $k|l$. Interchanging the roles of $k,l$ we get $l|k$ and so we must have $k=l$ since $k,l\geq 1$. 

Now consider $\sigma$-twisted conjugacy classes.  Since $A\sim_\sigma B$ if and only if $AJ=X(BJ)X^{-1}$ 
for some $X\in \SL(n,\bz)$. 
We need only show that the set $\{\tr(AJ)\mid A\in \SL(n,\bz)\}$ is infinite. 
Let $A'\in \SL(n-1,\bz)$ and let $A=diag(A',1)$ where $A'\in \SL(n-1,\bz)$.  Then $AJ=diag(A',-1)$.  Therefore 
$\tr(A)=\tr(A')-1$. Since $n>2$, the set $\{\tr(A')\mid A'\in \SL(n-1,\bz)\}$ is infinite and we conclude that  
$R(\sigma)=\infty$.

The proof of that $R(\sigma\tau)=\infty$ is similar and so we omit the details. 
\end{proof} 

It is possible to give a more direct proof of the $R_\infty$-property for $\SL(2,\bz)$ as for $\SL(n,\bz), n>2,$ given above, using the description of the (outer) automorphism group of $\SL(2,\bz)$ given in 
\cite[Theorem 2]{hr}.

It is known that the $R_\infty$-property is not inherited by finite index subgroups in general. For example, the infinite dihedral group, which contains the infinite cyclic group as an index $2$ subgroup, has the $R_\infty$ -property (whereas $R(-id_\bz)=2$). (See \cite{gw}.)  However we have the following result. 

Let $\Gamma_m$ denote the principal level $m$ congruence subgroup of $\SL(n,\bz)$; thus $\Gamma_m$ is the kernel of the surjection $\SL(n,\bz)\lr \SL(n,\bz/m\bz)$ induced by the surjection $\bz\lr \bz/m\bz$. 

\begin{theorem}\label{congruence} Let $n\geq 3$.   
Let $\Lambda$ be a non-central normal subgroup of 
$\SL(n,\bz)$.  Then $\Lambda$ has the $R_\infty$-property.   
\end{theorem}
\begin{proof}
Let $\Gamma=\SL(n,\bz)$.  We shall use the notations 
introduced in the proof of Theorem \ref{linear}. 
It is known that any finite index subgroup of $\SL(n,\bz)$ contains $\Gamma_m$ for some $m\geq 2$. 
This is the congruence subgroup property for $\SL(n,\bz), n>2$.  See \cite[\S4.4]{sury}.  

 Let $M=(m_{i,j})\in \SL(n,\bz)$ and let $\phi:=\phi_M$ be the restriction to $\Lambda$ of the inner automorphism $\iota_M$ of $\Gamma$.  Then $X\sim_\phi Y$ if and only if $XM=Z(YM)Z^{-1}$ for some $Z\in \Lambda$. 
In particular $\tr(XM)=\tr(YM)$.  To show that $R(\phi)=\infty$ we need only show that 
the set $\{\tr(AM)\mid A\in \Lambda\}$ is infinite for any $M\in SL(n,\bz)$.   
There are two cases to consider: (1) $m_{ii}\neq 0$ for 
some $i$, (2) $m_{ii}=0$ for all $i$.

\noindent 
{\it Case (1):}  Without loss of generality we may assume that 
$m_{11}\neq 0$.  
Let $k>1$ and let $X(k)$ be the 
block diagonal matrix $X(k)=diag(C(k), I_{n-2})$ where $C(k)=\left(\begin{array}{cc} k^2+1 & k\\k&1\\
 \end{array}\right)$. A straightforward computation shows that $\tr(X(k)M)=(k^2+1)m_{11}+k(m_{12} +m_{21})+\sum_{2\leq j\leq n}m_{jj}$.  Therefore 
 $\tr(X(k)M)=\tr(X(l)M)$ 
if and only if $(k+l)m_{11}+m_{12}+m_{21}=0$.  
Choose $k_0>(m_{12}+m_{21})/m_{11}$. Then  
$X(mk),k\geq k_0,$ belong to pairwise distinct $\phi$-twisted conjugacy classes in this case. 

\noindent
{\it Case (2):}  Without loss of generality assume that $m_{12}\neq 0$.   Let $A(k)$ be as in the proof of Theorem \ref{linear}. Then $\tr(A(k)M)=km_{1,2}$. Therefore 
$\tr(A(k)M)=\tr(A(l)M)$ if and only if $k=l$.  Since $A(mk)\in \Gamma_m\subset \Lambda$ for all $k$, it follows that $R(\phi)=\infty$ in this case as well. 
   
Suppose that $\tau(\Lambda)=\Lambda$ where 
$\tau(X)= {}^{t}X^{-1}$ as in the proof of Theorem \ref{linear}. 
We see that $R(\tau|\Lambda)=\infty$ 
arguing as we did to establish that $R(\tau)=\infty$ in 
the proof of Theorem \ref{linear} by considering the set of elements $A(mk)\in \Gamma_m\subset \Lambda, k\geq 1$. Similarly, we show that 
$R(\theta|\Lambda)=\infty$ for each representative $\theta$ of the outer automorphisms of $\Gamma$ which leaves $\Lambda$ invariant. 

To complete the proof, we need only show that {\it every} automorphism of $\Lambda$ 
extends to an automorphism of $\Gamma$.  
For this purpose we observe that 
the $\br$-rank of the semi simple Lie group $G:=\SL(n,\br)$ equals $n-1\geq 2$.  
Let $\theta:\Lambda\lr \Lambda$ be any automorphism. By the Mostow-Margulis strong rigidity theorem \cite[Chapter 5]{zimmer}, $\theta$ extends to an automorphism $\wt{\theta}:\SL(n,\br)\lr \SL(n,\br).$ 
By a result of Newman \cite[Lemma 2]{newman} we have $N_G(\Lambda)=\Gamma$. So 
$\wt{\theta}$ restricts to an automorphism $\bar{\theta}$ of $\Gamma$.  Thus $\theta$ is the restriction of an automorphism of $\Gamma$, namely $\bar{\theta}$. This 
completes the proof. 
\end{proof}

\begin{remark}{\em 
Recall that Fel'shtyn and Gon\c{c}alves \cite{fg} have shown that $\Sp(2n,\bz)$ has the $R_\infty$- property.  
One could also establish this result along the same lines as for $\SL(n,\bz)$ given above.  We assume that $n\geq 2$ as $\Sp(2,\bz)=\SL(2,\bz)$. 
To fix notations, 
regard $\Sp(2n,\bz)$ as the subgroup of $\SL(2n,\bz)$ 
which preserves the skew symmetric form $\beta:\bz^{2n}\times \bz^{2n}\lr \bz$ defined as $\beta(e_{2i},e_{2j})=0=\beta(e_{2i-1},e_{2j-1}),~
\beta(e_{2i-1},e_{2j})=\delta_{ij}, 1\leq i\leq j\leq n$ (Kronecker $\delta$).  Equivalently 
$\Sp(2n,\bz)=\{X\in \SL(2n,\bz)\mid {}^tXJ_0X=J_0\}$ 
where $J_0=diag(j_0,\ldots, j_0), j_0:=
\left(\begin{array}{cc} 0 & 1\\-1&0\\
 \end{array}\right)$ is the matrix of $\beta$. 
Infinitude of (untwisted) conjugacy classes follows from the residual finiteness of $\Sp(2n,\bz)$. (cf. Proposition \ref{residual}).  Alternatively, observe that $X(k)\in \Sp(2n,\bz)$ 
where $X(k)$ is as in the proof of Theorem \ref{congruence}.  This shows that the trace function is unbounded on $\Sp(2n,\bz)$. 

To complete the proof, we need only verify that $R(\phi)=\infty$ for representatives of the elements of $\out(\Sp(2n,\bz))$.     
One knows from \cite{reiner} that the outer automorphism 
group of $\Sp(2n,\bz)$ is isomorphic to 
$\bz/2\bz$ 
if $n> 2$ and is isomorphic to $\bz/2\bz\times \bz/2\bz$ 
when $n=2$. 

The generators of the outer automorphism groups may be 
described as follows. Let 
$\theta$ be the automorphism of $\Sp(2n,\bz)$ which is conjugation  by $J:=diag(I',I_{2n-2})\in \gl(2n,\bz)$ where $I'=
\left(\begin{array}{cc} 0 & 1\\1&0\\
 \end{array}\right).$ Let $\phi$ be the 
 automorphism of $\Sp(4,\bz)$ defined as $\phi(X)=\chi(X)X$ where $\chi:\Sp(4,\bz)\lr \{1,-1\}$ is the (non-trivial) central character.  Then $\out(\Sp(2n,\bz))=\langle\theta\rangle, n>2,$ and $\out(\Sp(4,\bz))=\langle\theta,\phi\rangle$.

To see that $R(\theta)=\infty$ we note that $\tr(X(k)J)=
2k+(2n-2)$.  Therefore the $X(k),~k\geq 1,$ belong 
to pairwise distinct $\theta$-twisted conjugacy classes.

As observed already in \cite[Lemma 3.1]{fg}, any $\phi$-twisted conjugacy class of $X$ is a union of 
the (untwisted) conjugacy class of $X$ and of $-X$.   
Since the number of conjugacy classes in $\Sp(4,\bz)$ is infinite, it follows that $R(\phi)=\infty$. 
Proof that $R(\theta\phi)=\infty$ is similar and omitted. This completes the proof. }
\end{remark}

It is an interesting problem to determine which 
(irreducible) lattices in semi simple Lie groups have 
the $R_\infty$-property.  
\section{Proof of the main theorem}
We now proceed to the proof of the main theorem.  Let $j:A\hookrightarrow \Lambda$ is the inclusion and $\eta:\Lambda\lr \Gamma$ the canonical quotient so that 
$1\lr A\hookrightarrow \Lambda\lr \Gamma\lr 1$ is 
an exact sequence of groups.

\noindent
{\it Proof of Theorem \ref{main}:}  
Let $\phi:\Lambda\lr \Lambda$ be any automorphism and let $f:A \lr \Gamma$ be the composition $\eta\circ\phi\circ j$. 
Note that since $A$ is normal in $\Lambda$, $\phi(A)$ is normal in $\Lambda$ and hence $f(A)$ is normal in $\Gamma$. 

(i)  In this case we claim that $f$ is trivial.  Suppose that 
$f(A)$ is {\it not} the trivial subgroup.  Since $\Gamma$ is a non-elementary group it does not contain a free abelian group of rank $2$. Since $\Gamma$ is torsion free, 
the centralizer of any non-trivial element of $\Gamma$ 
is infinite cyclic. By \cite[Corollary 3.10, Chapter III.$\Gamma$]{bh} $f(A)$ is quasi-convex. Hence by 
\cite[Proposition 3.16, Chapter III.$\Gamma$]{bh} the subgroup $f(A)$ has finite index in its 
normalizer, which is $\Gamma$. This contradicts the assumption that $\Gamma$ is non-elementary. 
Therefore $f(A)$ must be trivial.  This means that $\phi(A)\subset A$ and we have the following  
diagram in which the top horizontal sequence is exact:
\[
\begin{array}{rllll}
A&\hookrightarrow &\Lambda&\lr& \Gamma\\ 
\phi|A\downarrow &&\phi\downarrow&&\downarrow\bar{\phi}\\
A&\hookrightarrow &\Lambda&\lr& \Gamma.\\ 
\end{array}
\]
Now $\bar{\phi}$ is a surjection since $\eta\circ\phi$ is. 
Since $\Gamma$ is assumed to be torsion-free, by Sela's theorem \cite{sela}, $\Gamma$ is hopfian and so $\bar{\phi}$ is an isomorphism. Therefore $\phi(A)=A$. Hence $A$ is characteristic in $\Lambda$. 
Since $\Gamma$ has the $R_\infty$-property by \cite{ll} (cf. \cite{felshtyn1}), Lemma \ref{characteristic} now implies that $\Lambda$ has the $R_\infty$-property.  
 
(ii)  
The group $\Gamma$ is a lattice in one of the simple linear Lie groups $G=
\SL(n,\br), \pgl(n,\br),\\
 \Sp(2n,\br), \textrm{PSp}(2n,\br)$.  These Lie groups have centre 
a group of order at most $2$.    Also, $\Gamma$ is hopfian. 
First we consider the case $\Gamma=\SL(n,\bz), \psl(n,\bz), \pgl(n,\bz),n>2,$ or $\Sp(2n,\bz), 
\textrm{PSp}(2n,\bz), n>1,$ so that the corresponding Lie group $G$ has real rank at least $2$.  By the normal subgroup theorem of Margulis \cite[Chapter 8]{zimmer}, the subgroup $f(A)$ being normal in $\Gamma$ is either of finite index or  
is contained in the centre of $G$.  Since $A$ is abelian, $f(A)$ cannot be of finite index in $\Gamma$.  Hence $f(A)\subset Z(\Gamma)$ the centre of $\Gamma$ which is of order at most $2$.  First assume that $f(A)$ is 
trivial.   Then we have $\phi(A)\subset A$.  Using the fact that $\Gamma$ is 
hopfian, we conclude as above, that $A$ is characteristic. Now $\Gamma$ has the $R_\infty$-property 
by Theorem \ref{linear} in the case of $\SL(n,\bz), \psl(n,\bz), \pgl(n,\bz)$ and by the work of 
Fel'shtyn-Gon\c{c}alves \cite{fg} in the case of $\Sp(2n,\bz), \textrm{PSp}(2n,\bz)$ (cf. Lemma \ref{finitecharacteristic}(i)).  It follows as in case (i) that $\Lambda$ also has the $R_\infty$-property.  
Now assume that $f(A)=Z(\Gamma)\cong\bz/2\bz$.  Set $\bar{\Gamma}
=\Gamma/ Z(\Gamma)$ which is the lattice $\psl(n,\bz)$ or $\textrm{PSp}(2n,\bz)$ in 
the corresponding Lie group of adjoint type. Let $N=\eta^{-1}(Z(\Gamma))$. Clearly $N/A\cong Z(\Gamma)$. 
Now we have the exact sequence 
$N\stackrel{\wt{j}}{\hookrightarrow} \Lambda\stackrel{\bar{\eta}}{\lr}\bar{\Gamma}$ where $\bar{\eta}$ is the canonical quotient map. Now we claim that $N$ is characteristic. Indeed, 
let $\wt{f}:N\lr \bar{\Gamma}$ be defined as  $\bar{\eta}\circ\phi\circ\wt{j}$.  Again using Margulis' normal 
subgroup theorem, the fact that $N$ is virtually abelian forces $\wt{f}(N)$ to be contained in the centre of $\bar{\Gamma}$. Since $\bar{\Gamma}$ has 
trivial centre, we must have $\wt{f}(N)\subset N$. Now 
$\bar{\Gamma}$ is again hopfian (being finitely generated and linear).  As before, we conclude that 
$N$ is characteristic. By Lemma \ref{characteristic} applied to $\bar{\Gamma}$ we conclude that $\Lambda$ has the $R_\infty$-property. 

We now consider the case $\SL(2,\bz)\cong \Sp(2,\bz).$  Proceeding as above we see that $f(A)$ is a normal abelian subgroup of $\SL(2,\bz)$. We need only show that $f(A)\subset \{I,-I\}$.  Let $F\subset \SL(2,\bz)$ be a free group of finite index which is normal. Then $F\cap f(A)$ is trivial since any normal subgroup of $F$ is a non-abelian free group. Hence $f(A)$ is finite as it imbeds in 
the finite group $\SL(2,\bz)/F$. Let $C$ be the image of $f(A)$ in $\psl(2,\bz)\cong(\bz/2\bz)*(\bz/3\bz)$ under the natural quotient map.  Since any 
element of finite order is conjugate to the generator 
of $\psl(2,\bz)$ of order $2$ or that of order $3$ (see \cite[Theorem 2.7, Chapter IV]{ls}). Since 
$C$ is normal and finite, it follows easily that $C$ is trivial. Hence $f(A)\subset \{I,-I\}$.

(iii) 
By Theorem \ref{congruence} the $R_\infty$-property holds for $\Gamma$. The rest of the proof is as in case (ii) above and hence omitted.

(iv)
If $M$ is compact, then $\Gamma$ is 
a torsion-free hyperbolic group and our statement follows from part (i).  In 
any case, $\Gamma$ is a lattice in $G$, the group of orientation preserving 
isometries of the universal cover of $M$. Thus $G$ is a simple Lie group with trivial centre and real 
rank $1$.  In particular, $G$ is linear and so $\Gamma$ is residually finite.  Indeed $G$ is the identity component of the real points ${\bf G}_\br$ of the complex linear algebraic group ${\bf G}$ of adjoint type whose Lie algebra equals $\textrm{Lie}(G)\otimes_\br\bc$.  

If $M$ is non-compact then 
$\Gamma$ is relatively hyperbolic (with respect to the family of 
stabilizers of the cusps of $M$).  Fel'shtyn \cite[Theorem 3.3]{felshtyn2} has established the $R_\infty$- property for such groups $\Gamma$.  

Next we show that $f(A)$ is trivial.  
Let $Z\subset {\bf G}_\br$ be the Zariski closure of $f(A)$ and let $H$ be the normalizer of $Z$ in ${\bf G}_\br$.
Then $H$ is an algebraic subgroup which contains $\Gamma$. 
Since $\Gamma$ is Zariski dense in ${\bf G}_\br$ by the Borel density theorem \cite{raghunathan}, it follows that $H={\bf G}_\br$ and so $Z$ is normal in ${\bf G}_\br$. 
Since $Z$ is abelian and since $G$ is simple, it follows that $Z$ is finite and is contained in the centre of ${\bf G}_\br$.  Therefore $f(A)$ {\it equals} $Z\cap G$ and is contained in the centre of $G$. Since the centre of $G$ is trivial, we conclude that $f(A)=\{1\}$.   
The rest of the proof is as in the previous cases  
above.  \hfill $\square$

We conclude this paper with the following remarks.
\begin{remark}{\em 
(i) Theorem \ref{main} contains as special cases 
the direct product $A\times \Gamma$ as well as the 
the restricted wreath product $C\wr \Gamma=(\oplus _{\gamma\in\Gamma}C_\gamma)\ltimes \Gamma$ where 
$C_\gamma=C$ is any cyclic group.
(ii) Let $P$ be any set of primes containing $2$; thus any homomorphism $A(P)\lr \bz/2\bz$ is trivial.  Let $A(P)=\bz[1/p|p\in P]\subset \bq$.
Note that $i: A(P)\lr A(Q)$ is any non-trivial homomorphism, 
then $P\subset Q$.      
Set $\Lambda(P):=A(P)\wr \Gamma$ where $\Gamma$ is as in Theorem \ref{main}. Suppose that $\theta: \Lambda(P)\lr \Lambda(Q)$ is an isomorphism. Then, as in 
the proof of Theorem \ref{main}, the composition 
$\oplus_{\gamma\in \Gamma}A(P)\hookrightarrow \Lambda(P)\stackrel{\theta}{\lr} \Lambda(Q)\lr \Gamma$ is trivial.  It follows that $\theta(\oplus_{\gamma\in \Gamma}A(P))\subset \oplus_{\gamma\in \Gamma}A(Q)$ and so 
$P\subset Q$. Similarly $Q\subset P$ and so $P=Q$. 
It follows that there are $2^{\aleph_0}$ many pairwise non-isomorphic countable groups $\Lambda$ satisfying the $R_\infty$-property 
for each $\Gamma$ 
as in Theorem \ref{main}.  The same conclusion can also 
be arrived at by considering the groups $A(P)\times \Gamma$. \\
}
\end{remark}


\begin{thebibliography}{99}
\bibitem{brauer} R. Brauer,
Representations of finite groups. 1963 Lectures on Modern Mathematics, Vol. I 133--175 Wiley, New York. 
\bibitem{bh} M. Bridson and A. Haefliger, {\it Metric spaces of nonpositive curvature,} Grundlehren der Mathematischen Wissenschaften, {\bf 319} Springer-Verlag, Berlin, 1999. 
\bibitem{felshtyn1} A. Fel'shtyn, 
The Reidemeister number of any automorphism of a Gromov hyperbolic group is infinite. Zap. Nauchn. Sem. S.-Peterburg. Otdel. Mat. Inst. Steklov. (POMI) {\bf 279} (2001), Geom. i Topol. {\bf 6}, 229--240.  
\bibitem{felshtyn2} A. Fel'shtyn,  New directions in Nielsen--Reidemeister theory.  Topology Appl. {\bf 157} (2010) 1724--1735.
\bibitem{fg} A. Fel'shtyn and D. L. Gon\c{c}alves, Twisted conjugacy classes in symplectic groups, mapping class groups and braid groups. With an appendix written jointly with Francois Dahmani. Geom. Dedicata {\bf 146} (2010), 211--223.
\bibitem{fh} A. Fel'shtyn and Hill, 
The Reidemeister zeta function with applications to Nielsen theory and a connection with Reidemeister torsion, K-Theory {\bf 8} (1994) 367--393.
\bibitem{gw} D. L. Gon\c{c}alves and P. Wong, Twisted 
conjugacy classes in nilpotent groups. J. reine angew. Math. {\bf 633} (2009) 11--27.
\bibitem{hr}  L. K. Hua and I. Reiner, 
Automorphisms of the unimodular group. Trans. Amer. Math. Soc. {\bf 71} (1951) 331--348.
\bibitem{ll} G. Levitt and M. Lustig,  Most  automorphisms of a hyperbolic group have very simple dynamics. Ann. Sci. Ecol. Norm. Sup. {\bf 33} (2000) 507--517.
\bibitem{ls} R. Lyndon and P. Schupp,
{\it Combinatorial group theory}. Reprint of the 1977 edition. Classics in Mathematics. Springer-Verlag, Berlin, 2001.
\bibitem{omeara} O. T. O'Meara,The automorphisms of the linear groups over any integral domain. 
Jour. rein. angew. Math. {\bf 223} (1966) 56--100.  
\bibitem{newman} M. Newman, Normalizers of modular groups, Math. Ann. {\bf 238} (1978) 123--129.
\bibitem{osin} D. Osin, Small cancellations over relatively hyperbolic groups and embedding theorems. Ann. of Math. (2)  {\bf 172} (2010),  1--39,
\bibitem{pyber}  L. Pyber, Finite groups have many conjugacy classes. 
J. London Math. Soc. {\bf 46} (1992) 239--249. 
\bibitem{raghunathan} M. S. Raghunathan, {\it Discrete subgroups of Lie groups}, Ergebnisse der Mathematik und ihrer Grenzgebiete {\bf 68}, Springer-Verlag, Berlin, 1972.
\bibitem{reiner} I. Reiner,
Automorphisms of the symplectic modular group. Trans. Amer. Math. Soc. {\bf 80} (1955) 35--50. 
\bibitem{sela} Z. Sela, 
Endomorphisms of hyperbolic groups. I. The Hopf property. 
Topology {\bf 38} (1999) 301--321. 
\bibitem{serre2}J.-P. Serre, {\it Trees} 
Translated from the French original by John Stillwell. Corrected 2nd printing of the 1980 English translation. Springer Monographs in Mathematics. Springer-Verlag, Berlin, 2003.
\bibitem{sury} B. Sury, {\it Congruence subgroup 
property},  Texts and Readings in Mathematics, {\bf 24}. Hindustan Book Agency, New Delhi, 2003.
\bibitem{zimmer} R. J. Zimmer, {\it Ergodic theory and semisimple Lie groups}, Monographs in
Mathematics {\bf 81}, Birkh\"{a}user, 1984.
\end{thebibliography}
\end{document}